\documentclass [a4paper,twoside]{article}
\usepackage{latexsym}
\title{ Stochastic perturbations of iterations of a simple,
  non-expanding, nonperiodic, piecewise linear, interval-map} 
\author{ Thomas Kaijser
\\
\it\small Link\"{o}ping University,
 Sweden
\/{\rm ;} \it\small thomas.kaijser@liu.se }

\date{}
\begin{document}
\maketitle
\newtheorem{lem}{Lemma}[section]
\newtheorem{thm}{Theorem}[section]
\newtheorem{prop}{Proposition}[section]
\newtheorem{corr}{Corollary}[section]
\newtheorem{conjecture}{Conjecture}[section]
\newtheorem{definition}{Definition}[section]
\newtheorem{example}{Example}[section]
\newtheorem{condition}{Condition}[section]
\newtheorem{observation}{Observation}[section]

\begin{abstract}
Let \(g(x)=x/2 + 17/30 \pmod{1}\), let \(\xi_i, i= 1,2,...\) 
be a sequence of independent, identically distributed random variables 
with uniform distribution on the interval \([0,1/15]\), define
\(g_i(x)=g(x)+ \xi_i \pmod{1}\) and, for \(n=1,2,...\), define 
\(g^n(x)=g_n(g_{n-1}(...(g_1(x))...)).\) For \(x \in [0,1)\) 
let \(\mu_{n,x}\) denote the distribution of \(g^n(x)\). The purpose 
of this note is to show that there exists a unique probability measure
\(\mu\), such that, for all \(x \in [0,1),\)  
\(\mu_{n,x}\) tends to \(\mu\) as \(n \rightarrow \infty\). 
This contradicts a claim by Lasota and Mackey from 1987 stating 
that the process has an asymptotic three-periodicity.

\vspace{.5cm}

\vspace{.5cm}

\end{abstract}

\section{Introduction}\label{section_introduction}
Let  \(S=[0,1)\), let \(g: S \rightarrow S \) be defined by 
\begin{equation}\label{f}
g(x) = ax + b \;\; \pmod{1}
\end{equation}
where 
\begin{equation}\label{ab}
a= 1/2\;\; and \;\; b =17/30 .
\end{equation} 
Let \(\xi_n , \;n =1, 2, ... \) be a sequence of independent, identically 
distributed, random variables, define 
\(g_n :  S \rightarrow S\) by 
\[ g_n(x) = g(x) + \xi_n \; \pmod{1}\]
and define \(g^{(n)}: S \rightarrow S, \; n=1,2,..\) recursively by
\[
g^{(1)}(x) = g_1(x)
\]
\[
g^{(n+1)}(x) = g_{n+1}(g^{(n)}(x)), \; n=1,2,... .
\]
We write \(\xi^{(n)}=(\xi_1, \xi_2,...,\xi_n)\) and, if we want to
emphasize \(g^{(n)}(x)\)'s dependence of \(\xi_1, \xi_2,...,\xi_n\), we write
\[
g^{(n)}(x) = g^{(n)}(x; \xi^{(n)}).
\]
In the paper \cite{LM87} from 1987, A. Lasota and M. C. Mackey 
considered the process 
\(\{g^{(n)}(x), n=1,2,...\}\) for  two choices 
of the  sequence \(\{\xi_n, n=1,2,...\}\). 

The first case they considered was the case when  
\[Pr[\xi_n   = 0] =1 ,\; n=1,2,... .  \]
From a stochastic point of view this choice is somewhat 
artificial since  in this case
the sequence \(\{g^{(n)}(x), n=1,2...\}\) is a deterministic sequence. 
 Using  results from the paper \cite{Kee80} by J.P. Keener,  Lasota and Mackey  
concluded that when the parameters \(a\) and \(b\) in the expression
(\ref{ab}) are chosen such that \(a=1/2\) and  \(b=17/30\), then the sequence 
\(\{g^{(n)}(x), n=1,2,...\}\)
is a nonperiodic sequence for any  initial value \(x\). (For a more
explicit proof of this fact see \cite{Nak15}; especially page 465.)     

 Lasota and Mackey then also
considered the case when each of the stochastic variables 
\(\xi_n, n=1,2,...\)  has a uniform distribution on
the interval \([0,1/15]\). Using computer simulations they observed
that the distributions of the sequence 
\(g^{(n)}(\xi_0; \xi^{(n)}\)), where \(\xi_0\) has approximate uniform 
distribution on the interval [0,1), 
follow  a 3-periodic pattern already for \(n \geq 10 \).
(See \cite{LM87}, Figure 1   or \cite{LM94}, Figure 10.5.1.)

Thus, what Lasota and Mackey observed was  that, although  
a function is such that  it gives rise to a nonperiodic sequence of
numbers  when iterated, if - at each time epoch - 
the sequence of iterations is perturbed by a small stochastic
number,  then the distributions of the elements in the sequence may
show a  periodic pattern. They formulate this observation as follows:
\newline
\newline  
" ... . However, the surprising 
content of Theorem 1 ( of \cite{LM87}) 
is that even in a transformation S that has 
aperiodic limiting behaviour, 
the addition of noise will result in asymptotic periodicity.

This phenomenon is rather easy to illustrate  numerically by considering...".
(See \cite{LM87}, page 149.)
\newline
\newline

In the book \cite{LM94} from 1994 by Lasota and Mackey, 
the authors also present the example 
described above. Part of  the text in \cite{LM94} concerning this
example reads as follows: 
\newline
\newline
"Thus, in this example (the example above)  we have a noise induced
period three asymptotic periodicity". (See \cite{LM94}, section 10.5, page 323.)
\newline
\newline

This observed transition from an aperiodic behaviour to 
a periodic behaviour - thanks to stochastic perturbations -
 is certainly an interesting observation. 
However  this conclusion is not completely  true 
in the sense that in the long run the  3-periodicity will slowly 
disappear. What holds is that for any initial value \(x\) 
the distributions of the process
\(\{ g^{(n)}(x,\xi^{(n)}), n=1,2,...\}\) will tend to a unique limit
measure.   

\section{  Motivation}
Last year (2015),  an interesting paper by F. Nakamura called 
{\em Periodicity of non-expanding piecewise linear maps and 
effects of random noises} was published (see  \cite{Nak15}). 
Unfortunately though, in the last 
section of the paper, the author considers
the stochastic process described above and makes the same claim
as Lasota and Mackey. In fact, Nakamura even quotes the sentence from
 \cite{LM94}, that  was  mentioned above, verbatim.

 It thus seems that  still 29 years since 
the paper \cite{LM87} was published and 22 years since 
the book \cite{LM94} was published,  the fact 
that the claim made by Lasota and Mackey concerning the limit
behaviour
of the distributions of the stochastic process described above  
is not correct, has not been pointed out in
the literature. This is the motivation to write down a proof 
of the fact  that the stochastic process considered by Lasota and
Mackey in \cite{LM87}, section 5,  and in \cite{LM94} section 10.5, 
has a unique limit distribution.

The proof  presented  below
is in principal quite straightforward and not difficult, but writing 
down all the details requires a few pages.

At this point it is worth mentioning that although the convergence rate to 
the unique limit measure is exponential -  that is of order \(O(\rho^n)\)
 where \(\rho < 1 \) -   the parameter \(\rho\) is yet 
so close to unity that it is quite likely that it will not be possible to
reach the limit distribution - nor even come close to the limit  distribution - 
by computer simulations.

The observation made by Lasota and
Mackey, that stochastic perturbation may induce a high degree of 
periodicity may certainly - under certain 
circumstances -  be a useful and valuable observation. 

\section{Some simple formulas} For \(17/30 \leq b \leq 19/30 \) define 
\(g_b:[0,1) \rightarrow [0,1)\) 
\begin{equation}\label{gb}
g_b(x)= x/2 + b \;\pmod{1}.
\end{equation}
From (\ref{gb}) follows that 
\[
g_b(x) = x/2 + b, \;\; if\;\; 0 \leq x < 2(1-b)
\]
\[
g_b(x)= x/2 +b -1 \;\; if\;\; 2(1-b)\leq x < 1.
\]
Next define \(g_b^{(n)}\) recursively by
\(g_b^{(1)} = g_b\), \( g^{(n+1)}_b = g_b\circ g^{(n)}_b\).
By  simple calculations we find that \(g^{(2)}(x)\) satisfies
\[
g_b^{(2)}(x) = x/4 + 3b/2, \;\; if \;\;0 \leq x < 4- 6b,
\]
\[
g_b^{(2)}(x) = x/4 + 3b/2 -1/2, \;\; if \;\;4-6b \leq x < 2(1-b),
\]
\[
g^{(2)}(x) = x/4 + 3b/2 -1, \;\; if\;\;2(1-b) \leq x < 1,
\]
and we find that \(g_b^{(3)}(x)\) satisfies
\[
g^{(3)}_b(x) = x/8 + 7b/4, \;\; if \;\;0 \leq x \leq 8-14b \;\;\; and 
\;\;\;17/30 \leq b < 4/7,
\]
\[
g^{(3)}_b(x) = x/8 + 7b/4 -1, \;\; if \;\;\max\{0,  8-14b\} \leq x < 4-6b
\]
\[
g^{(3)}_b(x) = x/8 + 7b/4 -1/2, \;\; if \;\;4-6b \leq x < 2(1-b)
\]
and
\[
g^{(3)}_b(x) = x/8 + 7b/4 -1/4, \;\; if \;\;2(1-b)\leq x < 1.
\]
 Note that if \(b\geq 4/7\) then the set 
\(\{x: 0 \leq x < 8-14b\;\} = \emptyset \). 

Next set \(A=[17/30, 1)\) and let \(I_A: S \rightarrow \{0,1\}\) denote the 
indicator function of \(A\) . The rotation number \(rot_{g_b}(x)\) of
\(g_b\) can be defined
by 
\[
rot_{ g_b}(x) = \lim _ {N \rightarrow \infty} N^{-1}\sum_{n=1}^N I_A(g_b^{(n)}(x))
\]
(See \cite{Kee80}, Definition 1.1, page 590.)  Since 
\(g_b(0) > \lim_{x \rightarrow 1} g_b(x) \) it follows from Lemma 3.1 of
\cite{Kee80} that \(rot_{g_b}(x)\) exists and is  independent of x. 

\begin{prop} If  \(4/7 \leq  b \leq 19/30 \) then 
\[
rot_{g_b}(x) = 1/3
\]
whereas if \(17/30 \leq b <4/7 \) then 
\[
rot_{g_b}(x) < 1/3.
\]
\end{prop}

We shall not prove this proposition since it will not be used 
in our proof of Theorem \ref{Nagaev} below. Let us just make a few observations.
\newline
1) If \(b=4/7=120/210\) then \(g_b^{(3)}(0) = 0\).
\newline
2) If \(b = 19/30\) and \(x_0 = 13/(7\cdot 15)\) then 
\(0 < x_0 < 4- 6(19/30)\) and 
\[
g^{(3)}_b(x_0)=x_0.
\]
3) The ratio between the sets  \([17/30 , 4/7] (=[119/210, 120/210]) \)  and 
\([4/7, 19/30](=[120/210,133/210] \) is equal to 1/13.   

The first two observations indicate the truth of the proposition.  
The third observation, that the ratio between the sets \([17/30 , 4/7]\)  and 
\([4/7, 19/30]\) is equal to \(1/13\) and thus quite small, 
 explains why computer simulations show a 3-periodic pattern. 
On the other hand, since the rotation number \(rot_{g_b}(x) < 1/3 \) when 
\(17/30 < b < 4/7\) it is not surprising that in the long run the sequence 
\(\{ g^{(n)}(x, \xi^{(n)}), n=1,2,...\}\) as defined in Section 1, has
a unique limit measure independent of \(x\), as we claimed above. 

We shall end this section stating yet one more relation which gives some more 
information about the mapping \(g_b:[0,1) \rightarrow [0,1)\) when \(b=17/30\). 

For, suppose that \(x =26/30- \epsilon\) where say for simplicity
\(0<\epsilon<1/100.\) Then, by simple calculations, we find that 
\[
g^{(4)}_b(26/30 -\epsilon) = 26/30 - \epsilon /16,
\]
if \(b=17/30\)
from which we see that \(g^{(4n)}_b(26/30 -\epsilon) \rightarrow
26/30\) as \(n\rightarrow \infty\) if \(b=17/30\), 
from which we can conclude  that \(g_b \) is  very close to a
{\em 4-periodic
function} if \(b=17/30\). That \(g_b\) is not a 
{\em 4-periodic function} when \(b=17/30 \) is easy to check by
showing that the equation \(g^{(4)}_b(x)-x = 0 \) has no solutions
when \(b=17/30\).

\section{A limit result}\label{section_limitresult}
Let \(S=[0,1)\), let \(\delta: S\times S \rightarrow S\) be defined by
\[
\delta(x,y) = |x-y|
\]
and let  \({\cal B}\) be the Borel field on \(S\) determined by \(\delta\).
Further, as before let
 \(g:S \rightarrow S\) be defined by 
\[g(x)=ax +b \;\pmod{1}.\] 
where \(a=1/2\) and  \(b=17/30\).

Next  let \(\Omega = [0,2/30]\), 
let \({\cal A}\) be the Borel field on \(\Omega\).
Set \(\Omega^1 = \Omega\), \({\cal A}^1 = {\cal A}\) and for \(n=2,3,...\)
define \(\Omega^n\) and \({\cal A}^n\) recursively by
\[
\Omega^n=\Omega \times \Omega^{n-1}
\]
\[
{\cal A}^n={\cal A}^{n-1}\otimes {\cal A}.
\]
We denote a generic element in \(\Omega^n\) by 
\(\omega^n =(\omega_1, \omega_2, ...,\omega_n)\).

Next let \(\{ f^{(n)}:S \times \Omega^{n} \rightarrow S, \;n=1,2,...\}\)
be a sequence of functions defined recursively by
\begin{equation}\label{f1}
f^{(1)}(x,\omega) =g(x)+ \omega \pmod{1}.
\end{equation}
\begin{equation}\label{fn}
f^{(n+1)}(x, \omega^{n+1})= f^{(1)}(f^{(n)}(x, \omega^n), \omega_{n+1}).
\end{equation}

Let   \(\{\xi_n, n=1,2,...\}\) be a sequence of independent, identically 
distributed, random variables having uniform distribution on the 
interval \(\Omega\) and set \(\xi^{(n)}=(\xi_1,\xi_2,...,\xi_n)\).
We denote the distribution of \(\xi_n\) by \(\lambda\) and the
distribution of \(\xi^{(n)}\) by \(\lambda^n\).

For \(n=1,2,...\) define \(K^n:S \times {\cal B} \rightarrow [0,1]\)
by
\begin{equation}\label{Kn}
K^n(x, A) = Pr[ f^{(n)}(x, \xi^{(n)})\in A] = \int_A f^{(n)}(x, \omega^n)\lambda^n(d\omega^n).
\end{equation}
\begin{thm} \label{Nagaev}
There exists a constant \(C > 0\), a constant \(\rho < 1\) 
and a measure \(\mu\) such that for all \(x \in S\) and all 
\(A \in {\cal B}\) 
\begin{equation}
 |K^n(x,A) - \mu(A)| \leq C \rho^n.
\end{equation}
\end{thm}

The proof can be regarded as a "routine matter". Our proof is 
based on a simple {\em coupling device}.

\section{ An auxiliary limit theorem for Markov chains}

Let \((S,{\cal F}, \delta)\) be a compact metric space where 
\({\cal F}\) is the Borel field induced by the metric \(\delta\).
Let \(P:S \times {\cal F}\rightarrow [0,1]\) be a transition
probability function (tr.pr.f). Let \(P^n:S \times {\cal F}
\rightarrow [0,1]\) denote the {\em n-step tr.pr.f } induced by 
\(P:S \times {\cal F}\rightarrow [0,1]\). Let 
\({\cal P}(S, {\cal F})\) denote the set of probability measures on
 \((S,{\cal F}))\). If \( \mu, \nu \in {\cal P}(S, {\cal F})\)
we let \(||\mu -\nu||\) denote the total variance distance
between \(\mu\) and \(\nu\) defined as usual by
\[
||\mu - \nu || = \sup\{\mu(F)-\nu(F): F \in  {\cal F}\} +
\sup\{\nu(F) -\mu(F):F\in {\cal F}\}
\]
and we let 
\({\tilde {\cal P}}(S^2, {\cal F}^2, \mu, \nu)\), denote the set
of all {\bf couplings} of \(\mu \) and  \(\nu\); that is 
the set of all probability measures \({\tilde \mu}\) on 
\((S\times S, {\cal F}\otimes {\cal F})\) such that
\[
{\tilde \mu}(F \times S) = \mu(F),\; \forall F \in {\cal F}
\]
and
\[
{\tilde \mu}(S \times F) = \nu(F), \; \forall F \in {\cal F}.
\]
We say that a tr.pr.f 
\({\tilde P}: S^2 \times {\cal F}^2 \rightarrow [0,1]\) 
is a {\bf Markovian coupling} of 
\(P:S \times {\cal F} \rightarrow [0,1] \) if 
for each \(x,y \in S\), \({\tilde P}(x,y, \cdot)\) is
a coupling of \(P(x, \cdot)\) and \(P(y,\cdot)\).

\begin{definition}\label{overlappingdefinition} We say that 
\(P:S \times {\cal F}\rightarrow [0,1]\) has the 
{\bf overlapping property} if there exists a set \(S_0 \in {\cal F}\) 
such  that 
\newline
1) there exist an integer \(N\) and a number \(\alpha_1>0\) such that
\[ 
\inf_{x \in S} P^N(x,S_0)\geq \alpha_1
\]
2) there exist a number \( \alpha_2 >0\) and a Markovian coupling
\({\tilde P}_0: S^2 \times {\cal F}^2 \rightarrow [0,1]\) 
of \(P\) such that  if \(D=\{(x,y) \in S \times S: x=y\}\) then
\[
\inf\{ {\tilde P}_0((x,y), D): x,y \in S_0\} \geq \alpha_2
\]
If we want to emphasize the parameters involved in the definition 
of the overlapping property, we
say that 
\(P:S \times {\cal F}\rightarrow [0,1]\) has the 
overlapping property with basic set \(S_0\), basic integer \(N_0\),
basic  coupling 
\({\tilde P}_0: S^2 \times {\cal F}^2 \rightarrow [0,1]\)
and basic lower bounds \(\alpha_1\) and \(\alpha_2\).
\end{definition} 

The following limit result holds.
\begin{thm}\label{helptheorem} Let \((S, {\cal F}, \delta)\) be
a compact metric space.
 Suppose \(P:S \times {\cal F}\rightarrow [0,1]\) has the
overlapping property. Then there exists a constant \(C > 0 \), a
constant \(0<\rho<1\) and a probability measure 
\(\mu \in {\cal P}(S,{\cal F})\), such that
\begin{equation}\label{Pestimate}
\sup\{||P^n(x, \cdot) - P^n(y,\cdot)||: x,y \in S\} \leq C \rho^n, \; n=1,2...
\end{equation}
and
\begin{equation}\label{muestimate}
\sup\{||P^n(x, \cdot) - \mu||: x\in S\} \leq C \rho^n, \; n=1,2...\; .
\end{equation}
\end{thm}
This theorem is not difficult to prove but for sake of completeness we give
a proof in the appendix.

\begin{corr}\label{corraftertheorem}
In order to prove Theorem \ref{Nagaev} it suffices to prove that
the tr.pr.f \(K:[0,1)\times {\cal B} \rightarrow [0,1]\) 
defined by 
\begin{equation}\label{K}
K(x, F) = Pr[f^{(1)}(x,\xi) \in F],
\end {equation} 
where \(f^{(1)}\) is defined by (\ref{f1}) and \(\xi\) is uniformly
distributed
on \([0, 1/15]\), has 
the overlapping property.
\end{corr}
{\bf Proof}. In order to be able to use Theorem 4.1 we need to verify that 
\(K^n: S\times {\cal B} \rightarrow [0,1]\), as defined by (\ref{Kn}), 
is in fact the \(n-step\) tr.pr.f
induced by the  tr.pr.f \(K:S\times {\cal B} \rightarrow [0,1] \)  defined by (\ref{K}). But 
this follows easily from the definition of \(\{ f^{(n)}:S \times
\Omega^{n} \rightarrow S, \;\;n=1,2,...\}\).
(See (\ref{f1}) and (\ref{fn}).) \(\;\Box\)

\section{Determining a basic set.} In order to prove that
the tr.pr.f \(K:[0,1)\times {\cal B} \rightarrow [0,1]\) defined by 
(\ref{K}) has the
overlapping property we shall first prove the following proposition.
\begin{prop}\label{propositionD}
 Let 
\[S_0=[0,3/30]
\]
and 
\[
D=\{(x,y)\in S \times S  : x=y \}. 
\]
Let \(K:[0,1)\times {\cal B} \rightarrow [0,1]\) be defined
as in Corollary \ref{corraftertheorem}. Then we can find 
a Markovian coupling 
\({\tilde K}:S^2 \times {\cal F}^2 \rightarrow [0,1] \) such that
\[
\inf \{{\tilde K}((x,y), D): x,y \in S_0\} \geq 1/4
\]
\end{prop}
{\bf Proof}.  
We devide \(S\times S\) inte four disjoint sets as follows.
 \[S_1= \{(x,y) \in S\times S : 0 \leq (y-x)/2 \leq 2/30\},
\]
\[S_2=\{(x,y) \in S\times S : 2/30 < (y-x)/2\},
\]
\[S_3= \{(x,y) \in S\times S : 0 < (x-y)/2 \leq 2/30\}
\]
and
\[
S_4 =\{(x,y) \in S \times S : 2/30< x-y)/2\}.
\]
As before let \(g(x)=x/2 +17/30 \pmod{1}\).

Next we define \(h_1:S \times S \times \Omega \rightarrow S \) by
\newline
a)
\[
h_1(x,y,\omega) =g(x) + \omega  + (y-x)/2\; \pmod{1}\;\;
\]
 if 
\[
\; (x,y) \in S_1 \;\;and \;\;\omega + (y-x)/2 \leq 2/30,
\]
b)
\[
h_1(x,y,\omega) = g(x) +\omega +(y-x)/2 - 2/30\; \pmod{1}
\]
if  
\[
 (x,y) \in S_1 \;\;and \;\;\omega + (y-x)/2 >2/30,
\]
and c)
\[
h_1(x,y,\omega)=g(x) + \omega \pmod{1}
\]
  if 
\[
 (x,y) \in S_2\cup S_3 \cup S_4,
\]
and we define  
\(h_2:S\times S \times \Omega \rightarrow S\) by
\newline
a) 
\[
h_2(x,y,\omega) = g(y) +\omega  + (x-y)/2 \;\pmod{1}
\]
 if 
\[
(x,y) \in S_3 \;\;and\;\; \omega + (x-y)/2 \leq 2/30,
\]
b)
\[
h_2(x,y,\omega) = g(y) +\omega +(x-y)/2 - 2/30\; \pmod{1}
\]
 if  
\[
 (x,y) \in S_3 \;\;and \;\;\omega + (x-y)/2 > 2/30, 
\]
and finally c)
\[
h_2(x,y,\omega)= g(y) + \omega \pmod{1}
\]
 if  
\[
(x,y) \in S_1\cup S_2 \cup S_4.
\]

 We also define 
\({\tilde h}=({\tilde h}_1, {\tilde h}_2): S\times S \times \Omega \rightarrow S\times S \)
by 
\[{\tilde h}_1(x,y,\omega)= h_1(x,y,\omega)\]
and
\[{\tilde h}_2(x,y,\omega)= h_2(x,y,\omega),\]
and we define 
\({\tilde K}:S\times S \times {\cal B}\otimes {\cal B}\rightarrow [0,1]\) by
\begin{equation}\label{Ktildedefinition}
{\tilde K}(x,y, F)=\lambda\{\omega: {\tilde h}(x,y,\omega) \in F\}
\end{equation}
\begin{lem}\label{Ktildecoupling} 
The function \({\tilde K}:S\times S \times {\cal B}\otimes
 {\cal B}\rightarrow [0,1]\)
defined above has the following properties.
\newline
a) \({\tilde K}:S\times S \times {\cal B}\otimes {\cal B}\rightarrow [0,1]\) is a tr.pr.f,
\newline
b) \({\tilde K}:S\times S \times {\cal B}\otimes {\cal B}\rightarrow [0,1]\) is a Markovian coupling of the tr.pr.f 
\newline \(K:S\times {\cal B} \rightarrow [0,1] \) defined by (\ref{Kn}),
\newline
c) if \(x,y \in S_0\) and \(D=\{(x,y) \in S \times S: x=y\}\)
then 
\begin{equation}\label{Lemc}
{\tilde K}(x,y,D) \geq 1/4.
\end{equation}
\end{lem}
{\bf Proof}. That \({\tilde K}(x,y, \cdot) \rightarrow [0,1]\) is a probability measure 
 for every \((x,y) \in S \times S\) follows
easily from the definition of \({\tilde K}:S\times S \times {\cal B}\otimes {\cal B}\rightarrow [0,1]\).
(See (\ref{Ktildedefinition}.) That also \({\tilde K}(\cdot, \cdot, F) \rightarrow [0,1]\) 
is \({\cal B}\otimes {\cal B}- measurable \) if \(F = A \times B\), 
where \(A\) and \(B\) are intervals, follows easily from the definitions of 
\(h_1:S\times S \times \Omega \rightarrow S\) and 
\(h_2:S\times S \times \Omega \rightarrow S\), and since the set 
of all rectangular sets
\(A \times B\) is a base for \({\cal B}\otimes {\cal B}\), it follows that
 \({\tilde K}(\cdot, \cdot, F) \rightarrow S\times S\) is \({\cal B}\otimes {\cal B}- measurable \)
for every \(F  \in {\cal B}\otimes {\cal B}\). Thus 
\({\tilde K}:S\times S \times {\cal B}\otimes {\cal B}\rightarrow [0,1]\) is a tr.pr.f 
which proves part a) of the lemma.

Next let us consider \({\tilde K}(x,y, A \times S)\) for \(A \in {\cal B}\).
From the definition of  \({\tilde K}:S\times S \times {\cal F}\otimes {\cal F} \rightarrow [0,1]\) 
(see (\ref{Ktildedefinition})) it follows that 
\[
K(x,y,A \times S) = \lambda\{\omega: h_1(x,y, \omega) \in A \}.
\]
If \((x,y) \in S_2 \cup S_3 \cup S_4\), then  
 \(h_1(x,y,\omega) = g(x) +\omega \pmod{1}\) from which immediately follows 
that in this case \({\tilde K}(x,y, A \times S) = K(x,A)\).

We also have  to consider the case when \((x,y) \in S_1\). In this case
\[
h_1(x,y,\omega) = g(x) + (y-x)/2 + \omega \pmod{1}  
\]
if \((y-x)/2 + \omega \leq 2/30 \)
and 
\[
h_1(x,y,\omega) = g(x) +(y-x)/2 +\omega - 2/30  \pmod{1}
\]
if \((y-x)/2 + \omega > 2/30\).  Now, if \(A \in {\cal B}\), 
and for each \(z \in [0, 2/30]\) we define 
\(A_z= \{\omega \in \Omega: g(x)+z+\omega \in A \;and \;z+ \omega < 2/30\}
\cup \{\omega \in \Omega: 
g(x) + z +\omega -2/30 \in A \;and \; z+\omega - 2/30 \geq 0\}\)
it follows easily that 
\(\lambda(A) = \lambda(A_z)\) from which follows that 
 \({\tilde K}(x,y, A \times S) = K(x,A)\) also in this case.

That  \({\tilde K}(x,y, S \times A) = K(y,A), \;\forall A \in {\cal
  B}\)
 can be proved in a similar way. Thereby part b) of the lemma is proved.

It remains to prove part c). But, if \(x,y \in S_0\) then 
\[
|y-x|/2 \leq 1/20.
\]
Suppose first that \(x \leq y\). We then find that 
\[
{\tilde h}_2(x,y, \omega) = g(y)+ \omega = 
y/2 +17/30 + \omega.
\] 
We also find that if  also \(0 \leq (y-x)/2 + \omega \leq 2/30\) then 
\[
{\tilde h}_1(x,y,\omega)= g(x)+ (y-x)/2 + \omega = 
x/2 + 17/30 + (y-x)/2 + \omega  = 17/30 + y/2 + \omega = \]
\[
g(y) +\omega = {\tilde h}_2(x,y, \omega).
\]
Hence if \(x,y \in S_0\), \(x \leq y\) and
\((y-x)/2 + \omega \leq 2/30\) then 
\[
{\tilde h}_1(x,y,\omega)={\tilde h}_2(x,y,\omega).
\]  
But clearly, since \(0 \leq (y-x)/2 \leq 1/20 < 2/30 \)
\[
\lambda \{\omega: (y-x)/2 + \omega \leq 2/30\} =
15 ( 2/30 - (y-x)/2) \geq 15(2/30 - 1/20) = 
\]
\[15(4-3)/60=1/4,
\]
from which follows that (\ref{Lemc}) holds if \(x,y \in S_0\)
and \(0\leq x \leq y \).
That (\ref{Lemc})  holds also if \(x,y \in S_0\)
and \(0\leq y < x \) can be proved similarly.
Thereby also part c) of
Lemma \ref{Ktildecoupling} 
is proved and from Lemma \ref{Ktildecoupling} follows  Proposition \ref{propositionD}.
\(\Box\).

\section{Finding return times for elements in the  basic set}
In the previous section we verified one of the two hypotheses
 that
the tr.pr.f \(K:S \times {\cal B} \rightarrow [0,1]\) has to  fulfill
in order to have the overlapping property. 
(See Definition \ref{overlappingdefinition}.) 
It thus remains to verify that we can find an integer \(N\) and a
number \(\alpha\) such that 
\[
\inf_{x \in S} K^N(x,S_0) \geq \alpha,
\]
where thus \(S_0=[0,3/30]\).

As a first step we shall in this section prove the following 
proposition.
\begin{prop}\label{proposition37}As above,  for \(n=1,2,...,\) let 
 \(K^n:S \times {\cal B} \rightarrow [0,1]\) be defined by
(\ref{Kn}) and let \(S_0=[0,3/30]\). There exist constants \(\alpha_0 >0\) and \(\beta_0 > 0\) such that 
\begin{equation}\label{K3}
\inf_{x \in S_0} K^{3}(x,S_0) \geq \alpha_0
\end{equation}
and
\begin{equation}\label{K7}
\inf_{x \in S_0} K^{7}(x,S_0) \geq \beta_0.
\end{equation}
\end{prop}
{\bf Proof.} Let 
\[
T_0 = \{ x: 0\leq x \leq 1/45\}
\]
and 
\[\Omega^3_0=\{\omega^3=(\omega_1,\omega_2, \omega_3) \in \Omega^3 :
\omega_1/4 + \omega_2/2 + \omega_3 <1/180\}.
\]
As before, for \(n=1,2,..\), let \(f^{(n)}:S\times \Omega^n \rightarrow S\) 
be defined by (\ref{f1}) and (\ref{fn}). Then, by simple calculations, we find that
\[
119/120 <
f^{(3)}(x, \omega^3)=x/8 +119/120 + \omega_1/4 + \omega_2/2 + \omega_3 <1
\]
if \(x \in T_0\) and \(\omega^{3} \in \Omega^3_0\).
Hence, if we define \(T_1=\{x : 119/120 \leq x <1\}\) we find that if \(x \in T_0\) then
\[
Pr[f^{(3)}(x, \xi^{(3)}) \in T_1] \geq Pr[\xi^{(3)} \in \Omega^3_0] 
\]
and since 
\[
 Pr[\xi^{(3)} \in \Omega^3_0] = (4/3)\cdot 15^3 \cdot(1/180)^3 =(4/3)\cdot
 (1/12)^3 
= (1/6)^4 = 1/1296,
\]
we find that
\[
Pr[f^{(3)}(x, \xi^{(3)}) \in T_1] \geq 1/1296
\]
if \(x \in T_0\).
Furthermore, since
\[
f(x,\omega) \in S_0
\]
if \(x \in T_1\) and \(\omega \leq 1/30\),
we find 
that
\[
Pr[f^{(1)}(x, \xi_1) \in S_0] \geq 1/2\]
if \(x\in T_1\) 
and hence
\begin{equation}\label{return4}
Pr[f^{(4)}(x, \xi^{(4)}) \in S_0] \geq (1/1296)\cdot(1/2) = 1/2592
\end{equation}
if \(x \in T_0\).

Next,  let \(x \in S_0\) and define
\[\Omega^3_x=\{\omega^3=(\omega_1,\omega_2, \omega_3) \in \Omega^3 :
x/8 +119/120 + \omega_1/4 + \omega_2/2 + \omega_3 \geq 1\}\]
and define 
\[W^3_x = \{\omega^3=(\omega_1,\omega_2, \omega_3) \in \Omega^3_x :
x/8 -1/120 + \omega_1/4 + \omega_2/2 + \omega_3 <1/45\}.\]

Then, since 
\newline
1)
\[
f^{(3)}(x, \omega^3) = x/8 +119/120 + \omega_1/4 + \omega_2/2 + \omega_3 -1
\]
if \(x\in S_0\) and \(\omega^{3} \in \Omega^3_x\cap W^3_x\),
\newline
2)
 \[-1/120 \leq x/8 - 1/120 \leq 1/240\;\; if \;\;\;x \in S_0\]
and 3)
\[1/45 - 1/240 =13/720 > 12/720 = 1/60,\]
it is not difficult to convince oneself that 
\begin{equation}\label{return3a}
Pr[f^{(3)}(x,\xi^{(3)}) \in T_0] \geq (4/3)\cdot 15^3 \cdot (1/60)^3 = 1/48
\end{equation}
if thus \(x \in S_0\).

Furthermore, since by monotonicity,
\[
\sup_{x \in S_0} Pr[f^{(3)}(x, \xi^{(3)}) \in T_1] =
\sup_{x \in T_0} Pr[f^{(3)}(x, \xi^{(3)}3) \in T_1] =\]
\[
 Pr[f^{(3)}(0, \xi^{(3)}) \in T_1] =(4/3)15^3 (1/120)^3 = 1/384
\]
it follows that we must have
\[
Pr[f^{(3)}(x, \xi^{(3)}) \in S_0] \geq 383/384
\]
if \(x \in S_0\).

By combining (\ref{return3a}) and (\ref{return4}) and using the fact
that \[f^{(7)}(x, \xi^{(7)}) = f^{(3)}(f^{(4)}(x, \xi^{(4)}), 
(\xi_5,\xi_6,\xi_7))\] it follows from the Markov property that 
\[
Pr[f^{(7)}(x, \xi^{(7)}) \in S_0] \geq (1/48) (1/2592) =1/(3^5 \cdot 2^9)> 1/125000 =
8\cdot 10^{-6} \]
if \(x \in S_0\).
Hence, setting 
\begin{equation}\label{alpha}
\alpha_0 =383/384
\end{equation}
and 
\begin{equation}\label{beta0}
\beta_0 = 1/(3^5\cdot 2^9)
\end{equation}
we find that (\ref{K3}) and (\ref{K7}) hold and thereby 
Proposition \ref{proposition37} is proved. \(\Box\)
\begin{corr}\label{12} Let \(\alpha_0\) and \(\beta_0\) be defined (\ref{alpha})
and (\ref{beta0}) respectively, let \(K:S\times {\cal B} \rightarrow [0,1]\)
be defined as in Corollary \ref{corraftertheorem}. Then, if \(x \in S_0\)
\[
K^{12}(x,S_0) \geq \alpha_0^4
\]
\[
K^{13}(x,S_0) \geq \alpha_0^2 \beta_0
\]
\[
K^{14}(x,S_0) \geq \beta_0^2
\]
\end{corr}
{\bf Proof.} Follows from (\ref{K3}), (\ref{K7}) and the Markov
property. 
\(\Box\).
\begin{corr} 
Let \(\alpha_0\) and \(\beta_0\) be defined (\ref{alpha})
and (\ref{beta0}) respectively, let \(K:S\times {\cal B} \rightarrow [0,1]\)
be defined as in Corollary \ref{corraftertheorem}. Then for every \(n
\geq 12\) there exists a number \(\gamma_n>0\) such that if \(x \in S_0\)
\[
K^n(x, S_0) \geq \gamma_n. 
\]
\end{corr}
{\bf Proof}. Follows  from Corollary \ref{12}, 
(\ref{K3}) and the Markov property. \(\Box\)

\section{ First entrance time to the basic set}
In the previous section we showed that 
\[
Pr[ f^{(n)}(x,\xi^{(n)}) \in S_0] > 0
\] 
for all \(n \geq 12\) if \(x \in S_0\). 
In this section we shall investigate \(K^n(x, S_0)\) 
when \(x \not \in S_0\).

We have already proved that 
\begin{equation}\label{T1}
K(x,S_0) \geq 1/2 \;\;if \;\;x \in T_1
\end{equation}
where thus \(T_1=[119/120, 1)\).

Next set \(T_2= [24/30,119/120]\). Since
\(f(24/30,\omega)= 12/30 +17/30 +\omega -1 \) if \(\omega \geq 1/30\)
and \(f(24/30, 2/30) = 1/30 < 3/30\) we find that 
\(f(24/30, \omega)\in S_0\) if \(\omega > 1/30\)
and since \(f(119/120, \omega) \in S_0\) if \(0\leq \omega < 1/30+1/240\),
we can conclude easily that 
\begin{equation}\label{T2}
K(x,S_0) \geq 1/2
\;\;
if \;\;x \in T_2.
\end{equation}
Next set \(T_3=[12/30, 24/30]\).
It is easily seen that in this case
\[
K(x, T_2\cup T_1]\geq 1/2.
\]
Hence 
\begin{equation}\label{T3}
K^2(x,S_0)\geq 1/4, \;\;if \;\;x \in T_3.
\end{equation}

It remains to consider the interval \(T_4=[3/30, 12/30]\).
This time it is easily seen that 
\[
K(x,T_3) \geq 1/2,
\]
and consequently 
\begin{equation}\label{T4}
K^3(x,S_0)\geq 1/8, \;\; if \; x \in T_4
\end{equation}
Combining (\ref{T1}),(\ref{T2}), (\ref{T3}) and (\ref{T4}) with 
Corollary \ref{12}, we can conclude 
that 
\[
\inf_{x \in S}K^{15}(x, S_0) \geq (1/2) \beta_0^2 \approx
3\cdot10^{-11}, 
\]
where thus 
\(\beta_0 =3^{-5}\cdot2^{-9}\approx 8\cdot 10^{-6} \).
Thereby we have verified that \(K:S \times{\cal B}\rightarrow [0,1]\)
has the overlapping property and hence Theorem \ref{Nagaev} follows from
Theorem \ref{helptheorem}. \(\Box \)

\section{Appendix 1. Proof of Theorem \ref{helptheorem}}
The purpose of this appendix is to prove Theorem \ref{helptheorem}. 
For sake of convenience we repeat the formulation.
\newline
\newline
{\bf Theorem \ref{helptheorem}}.
{\em  Let \((S, {\cal F}, \delta)\) be
a compact metric space. Suppose \(P:S \times {\cal F}\rightarrow [0,1]\) has the
overlapping property. Then there exists a constant \(C > 0 \), a
constant \(0<\rho<1\) and a probability measure $\mu$ such that
\[
\sup\{||P^n(x, \cdot) - P^n(y,\cdot)||: x,y \in S\} \leq C \rho^n, \; n=1,2...,
\]
and
\[\sup\{||P^n(x, \cdot) - \mu||: x\in S\} \leq C \rho^n, \; n=1,2...\; .
\]
}
{\bf Proof}. Let \(B[S,{\cal F}]\) denote the bounded, real,
Borel-measurable functions on \((S,{\cal F})\).  For 
\(u \in B[S,{\cal F}]\) define 
\[
||u|| = \sup\{ |u(x)|: x\in S\}
\]
and 
\[osc(u)= \sup\{u(x)-u(y): x, y \in S\}.
\]
For \(u \in B[S,{\cal F}]\) and \(\mu \in {\cal P}(S, {\cal F})\)
we write 
\[
\int_Su(x)\mu(dx) = \langle u, \mu\rangle. 
\] 

Next, let \(\mu, \nu \in {\cal P}(S, {\cal F})\).
It is well-known that
\begin{equation}\label{normdefinition}
||\mu - \nu|| = 
\sup\{ \langle u, \mu\rangle - \langle u, \nu\rangle : u \in B[S, {\cal F}], 
||u|| \leq 1\}.
\end{equation}

Thus, what we need to prove is that there exists a constant \(C\)
and a number \(0 < \rho < 1 \), such that for \(x,y \in S\)
\[
\sup\{\langle u, P^n(x, \cdot)\rangle -  
\langle u, P^n(y, \cdot)\rangle : u \in B[S, {\cal F}], 
||u|| \leq 1\} < C \rho^n,  n=1,2,...\;  .
\]

We start our proof with the following lemma.
\begin{lem}\label{xlem} Let \(\mu, \nu \in {\cal P}(S, {\cal F})\) 
and suppose that there exists a coupling \({\tilde \mu}\) of
\(\mu\) and \(\nu\) such that  
\[
{\tilde \mu}(D)=\alpha > 0
\]
where as above \(D=\{(x,y) \in S\times S : x=y\}\).
Let \(u \in B[S, {\cal F}]\).
Then 
\[|\langle u, \mu\rangle - \langle u, \nu\rangle| \leq (1-\alpha)osc(u).
\]
\end{lem}
{\bf Proof}. Let us first point out that the diagonal set \(D\)
belongs to the \(\sigma-field\) \({\cal F}\otimes {\cal F}\) since
\((S, {\cal F}, \delta)\) is a compact metric space. Next let 
\(u \in B[S, {\cal F}]\). Then
\[
|\int_S u(x) \mu(dx) - \int_S u(x)\nu(dx)|=
|\int_{S\times S}( u(x)-u(y)) \mu(dx)\nu(dy)|=\]
\[
|\int_{S\times S} (u(x)-u(y)) {\tilde \mu}(dx,dy)|\leq
\]
\[
|\int_{(S\times S)\setminus D} (u(x)-u(y)) {\tilde \mu}(dx,dy)|+
|\int_{D}( u(x)-u(y)) {\tilde \mu}(dx,dy)|\leq (1-\alpha)osc(u).
\]
\begin{corr}\label{xcorr1} Let \(P:S\times {\cal F} \rightarrow [0,1]\) be the 
tr.pr.f of Theorem \ref{helptheorem}. Since \(P\)
has the overlapping property there exist a basic set 
\(S_0\), a basic Markovian coupling \({\tilde P}_0\)
and a constant \(\alpha_2 > 0\) such that 
\[
\inf\{ {\tilde P}_0(x,y, D): x,y \in S_0\} \geq \alpha_2
\]
Let \(x,y \in S_0\). Then 
\[
||P(x,\cdot)-P(y,\cdot)|| \leq 1-\alpha_2 .
\]
\end{corr}
{\bf Proof}. Follows from Lemma \ref{xlem} and (\ref{normdefinition}). \(\Box\)

\begin{corr}\label{xcorr2} Let \(P\),
\(S_0\), \({\tilde P}_0\)
and \(\alpha_2 > 0\) be as in Corollary \ref{xcorr1}, and let 
\(\mu,\nu \in {\cal P}(S, {\cal F})\) be such that 
\[
\mu(S_0) \geq \alpha
\]
and
\[
\nu(S_0) \geq \alpha.
\]
Define \(\mu_1 \in {\cal P}(S, {\cal F})\) by
\[
\mu_1(F) = \int_SP(x,F)\mu(dx), \;\;F \in {\cal F}
\]
and
\(\nu_1 \in {\cal P}(S, {\cal F})\) by
\[
\nu_1(F) = \int_SP(x,F)\nu(dx), \;\;F \in {\cal F}.
\]
Then 
\[
||\mu_1-\nu_1|| \leq 1-\alpha_2\cdot \alpha^2. 
 \]
\end{corr}
{\bf Proof}. Define \({\tilde \mu}_1 \in {\cal P}(S^2, {\cal F}^2)\)
by
\[
{\tilde \mu}_1(A) = 
\int_{S\times S}{\tilde P}_0(x,y, A)\mu(dx)\nu(dy).
\]
It is easily checked that \({\tilde \mu}_1\) is a
coupling of \(\mu_1\) and \(\nu_1\). Furthermore we find that
\[
{\tilde \mu}_1(D) = 
\int_{S\times S}{\tilde P}_0(x,y, D)\mu(dx)\nu(dy) =
\int_{S_0 \times S_0}{\tilde P}_0(x,y,
D)\mu(dx)\nu(dy) +\]
\[
\int_{(S \times S)\setminus (S_0 \times S_0)}{\tilde P}_0(x,y,
D)\mu(dx)\nu(dy) \geq
 \alpha^2 \alpha_2 + 0.\]
From Lemma \ref{xlem} now follows that
\[|\langle u, \mu_1\rangle - \langle u, \nu_1\rangle| \leq 
(1-\alpha^2\alpha_2)osc(u)
\]
if \(u \in B[S, {\cal F}]\), which implies that
\[
||\mu_1-\nu_1|| \leq (1-\alpha_2\cdot \alpha^2). \;\;\;\Box
\]
\begin{corr} \label{xcorr3}Let \(P: S \times {\cal F} \rightarrow
  [0,1]\) have the 
{\em overlapping property} with basic set \(S_0\), basic integer \(N_0\),
basic  coupling 
\({\tilde P}_0: S^2 \times {\cal F}^2 \rightarrow [0,1]\)
and basic lower bounds \(\alpha_1\) and \(\alpha_2\). Then 
\[
||P^{N_0+1}(x, \cdot) - P^{N_0 +1}(y,\cdot)|| \leq (1- \alpha_1^2 \alpha_1),\;
\;\forall x,y \in S.
\]
\end{corr}
{\bf Proof.} Let \(x,y \in S\). Since 
\[
P^{N_0}(z, S_0) \geq \alpha_2, \;\forall z \in S
\]
it is clear that 
\[
P^{N_0}(x, S_0) \geq \alpha_2, 
\]
and 
\[
P^{N_0}(y, S_0) \geq \alpha_2.
\]
Since 
\(P^{N_0+1}(x, \cdot) \in {\cal P}(S, {\cal F})\) is defined
by
\[
P^{N_0+1}(x, F) = \int_SP(z,F)P^{N_0}(x,dz)
\]
and similarly 
\(P^{N_0+1}(y, \cdot) \in {\cal P}(S, {\cal F})\) is defined
by
\[
P^{N_0+1}(, F) = \int_SP(z,F)P^{N_0}(y,dz)
\]
we  see that the hypotheses of Corollary \ref{xcorr2} are satisfied.
The conclusion of Corollary \ref{xcorr3} now follows from Corollary \ref{xcorr2}. \(\Box\)

Next, let \(T:B[S, {\cal F}] \rightarrow B[S, {\cal F}]\) be defined by
\[
Tu(x)=\int_Su(y)P(x,dy)
\]
where thus \(P\) has the properties of the theorem under consideration.
If \(u \in B[S, {\cal F}]\) we may write 
\[
T^mu = u_m
\]
if convenient.

Next set\(N_1=N_0+1\) and \(\rho_1 = 1 - \alpha_1^2 \alpha_2\). 
From Corollary \ref{xcorr3}  it follows that 
\[
\sup \{T^{N_1}u(x) - T^{N_1} u(y): x,y \in S\} \leq \rho_1 osc(u)
\]
for all \(u \in B[S, {\cal F}]\).  Hence, for \(m=1,2,...\) 
\[
osc(T^{N_1+m}) \leq \rho_1 osc(u_m).
\]
By induction it follows that
\[
osc(T^{kN_1}) \leq osc(u) \rho_1^k, \; k=1,2...
\]
Since also \(osc(T^u) \leq osc(u), \forall u \in B[S, {\cal F}]\) we conclude
that 
\begin{equation}\label{oscestimate}
osc(T^nu) \leq C \rho^n osc(u), \;
\end{equation}
for all \(u \in B[S, {\cal F}]\)
if \(\rho\) and \(C\) are defined by 
\[
\rho = (\rho_1) ^{1/N_1}
\]
\[
C=1/\rho,
\] 
and since (\ref{oscestimate}) holds for all \(u
 \in B[S, {\cal F}]\), the estimate (\ref{Pestimate}) also holds and thereby
the first conclusion of Theorem \ref{helptheorem} is proved.
(See (\ref{Pestimate}).) 
 
That also the second inequality of Theorem \ref{helptheorem} holds,
follows easily from the first as follows. First, since 
\(osc(T^n(u)) \rightarrow 0\) and \((S, {\cal F}, \delta)\) 
is supposed to be a compact metric space, 
it follows that there exists a unique, invariant measure \(\mu\), such that 
\[
\lim_{n\rightarrow \infty} \int_Su(y)P^n(x,dy) - 
\langle u, \mu\rangle = 0, \;\; \forall x \in S.
\]
Furthermore, if we define 
\(Q:{\cal P}(S, {\cal F}) \rightarrow  {\cal P}(S, {\cal F})\) by
\[
Q\nu(A) = \int_S P(x, A) \nu(dx), \;\; \forall A \in {\cal F}
\]
and use the fact that if \(u \in B[K, {\cal F}]\) and  
\(\nu \in {\cal P}(K, {\cal F}\) then
\[
\langle Tu, \nu \rangle = \langle u , Q\nu\rangle 
\]
and the fact that \(\mu = Q\mu\) since \(\mu\) is invariant, 
we find that if \(u \in B[S, {\cal F}] \) then for \(n=1,2,...\) we have 
\[
|\int_S u(y)P^n(x,dy)- \langle u, \mu\rangle| = 
\]
\[
|\int_S (T^nu(x)-T^nu(y)) \mu(dy)| 
\]
which together with (\ref{oscestimate}) implies that for all 
\(x \in S\) 
\[
|\int_S u(y)P^n(x,dy)- \langle u, \mu\rangle|  
\leq osc(u) C\rho^n 
\]
which implies that
\[\sup\{||P^n(x, \cdot) - \mu||: x\in S\} \leq C \rho^n, \; n=1,2... 
\]
and thereby also the second 
 conclusion of Theorem \ref{helptheorem} is proved.
(See (\ref{muestimate}).) 
\(\;\;\Box\)

\end{document}